\documentclass{amsproc}
\usepackage{amssymb, latexsym}
\newcommand{\half}{\frac{1}{2}}

\newcommand{\e}{\epsilon}

\newcommand{\lam}{\lambda}
\newcommand{\Lam}{\Lambda}

\newcommand{\R}{{\mathbb R}}
\newcommand{\Rp}{{\mathbb R}^+}
\newcommand{\C}{{\mathbb C}}
\newcommand{\T}{{\mathbb T}}

\def\cal{\mathcal}

\newtheorem{theorem}{Theorem}[section]
\theoremstyle{definition}
\newtheorem{definition}[theorem]{Definition}
\theoremstyle{remark}
\newtheorem{remark}[theorem]{Remark}
\theoremstyle{proposition}
\newtheorem{proposition}[theorem]{Proposition}
\theoremstyle{lemma}

\theoremstyle{corollary}

\numberwithin{equation}{section}

\begin{document}

\title{KdV and Almost Conservation Laws}

\author{G.~Staffilani}
\address{Brown University and Stanford University}
\thanks{G.S. is supported in part by N.S.F. Grant DMS 0100375 
and a grant from Hewlett and Packard Foundation.}
\email{gigliola@math.brown.edu}


\begin{abstract}

In this article we  illustrate a new method  to extend local 
well-posedness results for dispersive equations to global ones. The 
main ingredient of this method is the definition of a family of what we call 
almost conservation laws. 
In particular we analyze the Korteweg-de Vries initial value problem 
and we illustrate in general terms how the ``algorithm'' that 
we use to formally generate almost conservation laws can be used to recover 
the infinitely many conserved integrals that make the KdV an integrable 
system.

\end{abstract}

\maketitle

\section{Introduction}
This short survey paper is concerned with a new method to prove global 
well-posedness results for dispersive equations below energy spaces, 
namely $H^{1}$ for the Schr\"odinger equation  and 
$L^{2}$ for the KdV equation. 

Even though I am the single writer of this article, all the new 
statements that I will make below have been proved together with
my collaborators J. Colliander, M. Keel, H. Takaoka and T. Tao. What 
started as a simple lunch at Stanford two years ago, evolved into a 
very fruitful collaboration in mathematics and a pleasant friendship. 
For whatever the reader appreciates in what follows, we all take the 
credit, for the mistakes, inaccuracies and the typos, I am the only one 
to blame!

Before starting with the story that I am set to tell, I
should warn the 
reader that because this article is a written version of the talk that 
I gave at the conference on Harmonic Analysis  in Mt. Holyoke College,
I will not present the complete proofs of the statements, but rather 
the main ideas involved in them.
The interested reader can check the references that I will list for a 
detailed proof of all the claims made. I also apologize in advance for not citing all the work that has been 
published in the context of well-posedness for dispersive equations. 
Here I will limit the bibliography to those publications that are in 
direct contact with the methods and the findings that I am about to 
describe.

We end this section with some notations.
Throughout the paper
we use $C$ to denote various constants. If $C$
depends on other quantities as well, this will be indicated by
explicit subscripting, e.g. $C_{\|u_0\|_2}$ will depend on 
$\|u_0\|_2$.  We use $A \lesssim B$ to denote an estimate of the
form $A \leq C B$, where $C$ is an absolute constant.  We use $a+$ and $a-$ to denote expressions of the
form $a+\varepsilon$ and $a-\varepsilon$, for some 
$0 < \varepsilon \ll 1$.

We use $\|f\|_{L^{p}}$ to denote the $L^{p}(\R)$ norm. For a fixed 
interval of time $[0,T]$ and a Banach space of functions $X$, we denote 
with $C([0,T],X)$ the space of the continuous maps from $[0,T]$ to $X$.

We define the spatial Fourier transform of $f(x)$ by
$$ {\cal F}(f)(\xi):= \hat f(\xi)
:= \int_{\R} e^{-i x \xi} f(x)\ dx$$
and the spacetime Fourier transform of $u(t,x)$ by
$$ {\cal F}(u)(\tau, \xi) :=
\hat u(\tau, \xi) := \int_{\R}\int_{\R} e^{-i (x \xi + t \tau)}
u(t,x)\ dt dx.$$
Note that the derivative $\partial_x$ is conjugated to multiplication
by $i\xi$ by the Fourier transform. We shall also define $D_x$ to be 
the operator conjugate to multiplication by 
$\langle \xi \rangle := 1 + |\xi|$.  We can then define the Sobolev
norms $H^s$ by
$$ \| f \|_{H^s} := \| D_x^s f\|_2 = \| \langle \xi \rangle^s \hat f
\|_{L^2_\xi}.$$

\section{Well-posedness and conservation laws}\label{con}

We consider the initial value problem (IVP) given by 
\begin{equation}\label{givp}
\left\{ \begin{array}{l}
\partial_tu+P(D)u+N(u)=0\\
  u(x,0) = u_{0}(x),
\end{array}\right.
\end{equation}
where $t\in \R$ and $x \in \R^{n}$ or $\T^{n}$, $P(D)$ is a 
differential operator with constant coefficients and $N(u)$ is the 
nonlinear part of the equation. For the moment we do not assume any 
special structure either for $P(D)$ or $N(u)$, we only assume that 
in terms of derivatives $P(D)$ is of at least one order higher than 
$N(u)$, in other wards we assume that the first equation in 
\eqref{givp} is {\em semilinear}. The function $u_{0}$ is called the 
initial profile and in general we assume that $u_{0}\in H^{s}$.

We will use the following definition for well-posedness:
\begin{definition}
The IVP \eqref{givp} is {\em locally well-posed} (l.w.p.) in $H^{s}$ if for any 
$u_{0}\in H^{s}$ there exists $T=T(\|u_{0}\|_{H^{s}})$ and a unique 
solution $u\in C([0,T],H^{s})$ for \eqref{givp}. Moreover the map that 
associates to each initial data its evolution is continuous.

We say that the IVP is {\em globally well-posed } (g.w.p.) in $H^{s}$ 
if for any $T>0$ the definition above is satisfied.
    \end{definition}
The question  of l.w.p is certainly the first one that one
investigates. After a positive result, then one trys to extend the 
local result to a global one. 

To convince the reader that proving well-posedness for a small interval 
of time is simpler than proving it for any fixed interval of large 
size, we briefly recall the contraction method.
We first use the Duhamel principle to write \eqref{givp} as the 
integral equation:
\begin{equation}\label{inteq}
    u(t,x)=W(t)u_{0}+\int_{0}^{t}W(t-t')N(u(t'))\,dt',
    \end{equation}
where $W(t)u_{0}(x)$ is the solution of the linear problem
\begin{equation}\label{lgivp}
\left\{ \begin{array}{l}
\partial_tv+P(D)v=0\\
  v(x,0) = u_{0}(x).
\end{array}\right.
\end{equation}
If one is willing to reduce the size of the interval of existence of 
the solution then one can replace \eqref{inteq} with
\begin{equation}\label{inteqd}
u(t,x)=\psi(t/\delta)W(t)u_{0}+
    \psi(t/\delta) \int_{0}^{t}W(t-t')N(u(t'))\,dt',
    \end{equation}
where $\psi(t)$ is a smooth cut-off function for the interval 
$[-2,2]$. We can still claim that $u$ solves \eqref{givp} in $[0,\delta]$ 
if and only if $u$ solves \eqref{inteqd} in the same interval. 
Now, consider the operator 
\begin{equation}\label{operator}
Lv(t,x)=\psi(t/\delta)W(t)u_{0}+
    \psi(t/\delta) \int_{0}^{t}W(t-t')N(v(t'))\,dt',
    \end{equation}
and assume that we are able to prove that there exists a Banach space $X^{s}$ 
and $s_{0}\in \R,$ such that for any $s\geq s_{0}$ we have 
$X^{s}\subset C(\R,H^{s})$  and 
\begin{eqnarray}
    \label{linear}& &\|\psi(t/\delta)W(t)u_{0}\|_{X^{s}}\leq 
    C_{0}\|u_{0}\|_{H^{s}},\\
\label{nonlin}  & &\left\|\psi(t/\delta) 
\int_{0}^{t}W(t-t')N(v(t'))\,dt'\right\|_{X^{s}}\leq 
C_{1}\delta^{\alpha}\|u\|_{X^{s}}F(\|u\|_{X^{s_{0}}}),\\
\label{nonlin12}& &\left\|\psi(t/\delta) 
\int_{0}^{t}W(t-t')[N(v_{1}(t'))-N(v_{2}(t'))]\,dt'
\right\|_{X^{s}}\leq\\  
\nonumber  && C_{1}\delta^{\alpha}
\max[\tilde F(\|v_{1}\|_{X^{s_{0}}}), \tilde F(\|v_{2}\|_{X^{s_{0}}})]
\|v_{1}-v_{2}\|_{X^{s}},
\end{eqnarray}
where $\alpha>0$, and $F, \tilde F:\R\longrightarrow \Rp$ are  
functions bounded on bounded sets.
If we set $a=2C_{0}\|u_{0}\|_{H^{s_{0}}}$ and we take 
$\delta^{\alpha}=1/4(C_{1}\max[ F(a), \tilde F(a)])^{-1}$, then 
the operator $L$
defined in \eqref{operator}
maps the ball $B_{a}$ in $X^{s}$ centered at the origin and radius $a$
into itself and is a contraction. Hence a unique fixed point exists 
and this is the unique solution for \eqref{givp}. Using a combination 
of \eqref{linear} and \eqref{nonlin12} one also obtains, for free, the 
continuity with respect to the initial data.
 
Arguably, this method has been used to prove the best 
 results on local well-posedness for a variety of dispersive 
equations (see \cite{B1}, \cite{B2}, \cite{KPV2}, \cite{KPVBilin}, 
and \cite{BouCol97}, just to name a few). 

We assume now for simplicity that $N(u)$,
the nonlinear part of the equation, is polynomial and that again
\eqref{linear}, \eqref{nonlin} and \eqref{nonlin12} still hold. Then 
the method we just described gives
well-posedness in $H^{s}, \, s\geq s_{0}$ in an interval of time $[0,T]$ such that
\begin{equation}\label{tu0}
T=C\|u_{0}\|_{H^{s_{0}}}^{-\beta},
\end{equation}
for some $\beta>0$.  We discuss  now
how to extend this short time result to a long time one. 

Accordingly to \eqref{tu0},  if we are willing to restrict our result to 
data small in $H^{s}$, then we can enlarge the time of existence. But 
this is not our goal here. We are looking in fact for a long time 
well-posedness for any initial data in $H^{s}$!. 

The first attempt that one can try is to iterate the short time result. 
Again by looking at \eqref{tu0}, it is clear that the obstacle in doing 
so will be the growth of $H^{s_{0}}$ norm of the solution $u(t)$ of 
\eqref{givp}. It is 
at this stage that uniform bounds for the Sobolev norms of the 
solution $u$ are needed and the conservation laws are the first 
source for such bounds. 

The existence of useful conservation laws depends on the structure of 
the equation in \eqref{givp}. So to continue our general exposition in 
this first section  we do not  write explicitly any conservation laws
involving the solution $u$, 
but instead we assume a consequence of them, whenever they are 
available, that is we assume that there exists $s_{*}\in \R$ such that
\begin{equation}\label{unifbound}
    \|u(t)\|_{H^{s_{*}}}\leq C^{*},
    \end{equation}
where $C^{*}$ does not depend on $t$. If  now $s_{*}\geq s_{0}$, then by 
\eqref{tu0}  and \eqref{unifbound} we can take 
$T^{*}=C(C^{*})^{-\beta}$ and iterate the local well-posedness result 
presented above. In the rest of the paper we will refer to this as 
{\em the method of conservation laws}.

We consider now two special examples of the IVP \eqref{givp}.
We start with the cubic defocusing Schr\"odinger equation in $\R^{2}$:
\begin{equation}\label{shr}
\left\{ \begin{array}{l}
i\partial_tu+\Delta u-|u|^{2}u=0,\\
  u(x,0) = u_{0}(x).
\end{array}\right.
\end{equation}
There are two conservation laws for this problem: the Hamiltonian
\begin{equation}
    \label{ham} \int_{\R^{2}}\frac{1}{2}|\nabla u|^{2}(x,t)
    +\frac{1}{4}|u|^{4}(x,t)\,dx=C_{1},
\end{equation}
and the $L^{2}$ norm
    \begin{equation}
    \label{l2}\int_{\R^{2}}|u|^{2}(x,t)\,dx=C_{0}.
\end{equation}
Using the Gagliardo-Nirenberg inequality, \eqref{ham} and \eqref{l2}, 
one obtains \eqref{unifbound} for $s_{*}=1$. On the other hand one can 
prove\footnote{and this is sharp!} 
that if $s>0$ then the IVP \eqref{shr} is well-posed in $H^{s}$ 
for an interval of time $[0,T]$, where $T\lesssim 
\|u_{0}\|_{H^{s}}^{-\beta}$, for some $\beta>0$, 
(see \cite{CW} and \cite{B1}). Then by the method of conservation laws
presented above one obtains global well-posedness for $s\geq 1$. 
So this method leaves the gap $s\in (0,1)$ open because the l.w.p., 
in the sense defined here, is barely missed at $s=0$, where the next 
conservation laws \eqref{l2} could have been used!

Next we pass to the KdV initial value problem 

\begin{equation}\label{kdv}
\left\{ \begin{array}{l}
\partial_tu+\partial_{x} u+\frac{1}{2}\partial_{x}u^{2}=0,\\
  u(x,0) = u_{0}(x),
\end{array}\right.
\end{equation}
where $x\in \R$ or $\T$. The KdV equation is special, in fact it 
enjoys infinitely many conserved integrals. Here we recall only the 
first four of them (notice that here $u$ is a real function!):
\begin{eqnarray}
    & &\int u(x,t)\,dx=C_{0}\\
\label{l2kdv}    & &\int u^{2}(x,t)\,dx=C_{1}\\
 \label{hamkdv}  & & \int \partial_{x}u^{2}+\frac{2}{3}u^{3}
 =C_{2}\\
 & & \int (\partial_{x}^{2}u)^{2}-\frac{5}{3}u(\partial_{x}u)^{2}
 +\frac{5}{9}u^{4}=C_{3}.
   \end{eqnarray}
Bourgain proved local well-posedness in $L^{2}$, hence by 
\eqref{l2kdv} and the method of conservation laws, 
global well-posedness in $L^{2}$ \cite{B1}. Later Kenig, Ponce 
and Vega showed that the IVP \eqref{kdv} on the line is locally 
well-posed in $H^{s}, s>-3/4$, leaving the gap $s\in (-3/4,0)$ open for 
global well-posedness. Similarly they proved that KdV
on the circle is locally well-posed in $H^{s}, \, s\geq -1/2$, 
leaving here the gap $s\in [-1/2,0)$. 

Similar results to the ones presented here for 
the IVP \eqref{shr} and \eqref{kdv} are available,
 with the obvious changes, also  for the modified KdV equations
\cite{B1} \cite{KPV2},
 the 1D Schr\"odinger equation with derivative nonlinearity
\cite{ozawa}, the KP-II equation \cite{B} and the 
Zakharov system \cite{BouCol97}.

\begin{remark}\label{remark1}
The method of conservation laws has two types of limitations. 
In general they only provide bounds for the $H^{1}$ norm (coming from 
the Hamiltonian), and the $L^{2}$ norm. Hence
when good local results are available  for rough data\footnote{
So far this has only been proved in low dimensions.}, 
these uniform bounds 
are not enough to cover all the possible indices $s$, and gaps 
are left as we showed above. The second limitation is that 
in higher dimensions well-posedness results are available only 
for relatively smooth data (in general in $H^{s_{n}}$, where $n$ is
the dimension and $s_{n}>n/2$). Then again uniform bounds in $H^{1}$
and $L^{2}$ are not enough (at least not yet!) to control these higher 
Sobolev norms. 
\end{remark}

\section{The method of Bourgain}
The method that we are about to describe is used to prove global 
well-posedness for rough initial data in low dimensions.  
It partially solves the first limitation of the method of conservation 
laws discussed in Remark \ref{remark1}. This method was first 
introduced by Bourgain \cite{Brefine} who considered the cubic, 
defocusing NLS on $\R^{2}$, ( but soon the reader will appreciate its 
generality). \\
As recalled above, for this IVP,  the method of conservation laws
leaves the gap $(0,1)$ between l.w.p and g.w.p.. So assume that 
 $u_{0}\in H^{s}$ and $s<1$. We split $u_{0}=\phi_{0}+\psi_{0}$, 
such that
\begin{equation}\label{cut}
\widehat{\phi_{0}}(\xi)=\chi_{\{|\xi|\leq N\}}\widehat{u_{0}}(\xi)
\hspace{1cm}\widehat{\psi_{0}}(\xi)=\chi_{\{|\xi|> N\}}
\widehat{u_{0}}(\xi),
\end{equation}
that is we decompose $u_{0}$ into low and high frequency parts. 
One can immediately observe that the low frequency part $\phi_{0}$
is smoother, but has a large norm:
$$
\|\phi_{0}\|_{H^{1}}\lesssim N, $$
while the high frequency part $\psi_{0}$ clearly does not improve its 
smoothness, but its  lower order norms are small:
$$
\|\psi_{0}\|_{H^{\sigma}}\lesssim 
N^{\sigma-s}, \, \, \mbox{ for }\, \, \sigma\leq s.$$
Then we evolve these two initial data. We call $u^{0}$ 
the evolution of the
low frequency part $\phi_{0}$ under the equation in \eqref{shr}. 
We call $v^{0}$ the evolution of the
high frequency part $\psi_{0}$ under the difference equation
$$i\partial_{t}v^{0}+\Delta v^{0}=|v^{0}+u^{0}|^{2}(v^{0}+u^{0})
-|u^{0}|^{2}u^{0}.$$
We can rewrite $v^{0}(t,x)=
S(t)\psi_{0}(x)+w^{0}(t,x)$, where $e^{it\Delta}\psi_{0}(x)$ 
is the solution of the associated linear problem 
\begin{equation}\label{shr0}
\left\{ \begin{array}{l}
i\partial_tv+\Delta v=0,\\
  v(x,0) = \psi_{0}(x).
\end{array}\right.
\end{equation}
and 
$$w^{0}(t,x)=\int_{0}^{t}S(t-t')(|v^{0}+u^{0}|^{2}(v^{0}+u^{0})
-|u^{0}|^{2}u^{0})\,dt'$$
is the nonlinear part. Clearly
$u^{0}(t)+ v^{0}(t)=u(t)$, where $u$ is the solution of 
\eqref{shr}.
There are  two key parts in Bourgain's argument. The first is 
that there exists 
$\delta=\delta(\|\phi_{0}\|_{H^{1}})>0$ such that 
    both $u^{0}(t)$ and $v^{0}(t)$ are defined for $t\in [0,\delta]$.
The second, more surprising\footnote{The reader should appreciate the 
remarkable fact that the nonlinear 
part of the evolution of the high frequency of the initial 
data is smoother than 
the data itself and small in the energy norm, hence it can be treated 
as an error!}, is that 
$$\|w^{0}\|_{H^{1}}\lesssim \frac{1}{N^{\alpha}},\, \, \mbox{ 
    for some } \alpha=\alpha(s)>0.$$
This is now the right set up for iteration. At this point we know that 
the unique solution $u(x,t)=u^{0}(x,t)=v^{0}(x,t)$ lives for all times 
in $[0,\delta]$. To proceed from $\delta$ to $2\delta$ we start a new 
IVP at time $\delta$ by assigning the new initial data
\begin{eqnarray*}
    \phi_{1}&=&u^{0}(\delta)+w^{0}(\delta)\\
    \psi_{1}&=&e^{i\delta\Delta}\psi_{0}
\end{eqnarray*}
and we repeat the argument above. An iteration like this would work 
on any finite interval $[0,T]$, as long as the total error is at most 
comparable with the size of $\|\phi_{0}\|_{H^{1}}$, the quantity that 
defines $\delta$, that is 
\[\sum_{1}^{M}\|w^{i}\|_{H^{1}}\sim \|\phi_{0}\|_{H^{1}}\sim N,\]
where $M\sim \delta^{-1}T$. By simple calculations on the 
explicit formula for $\delta$ and $\alpha(s)$ that we do not report 
here, one obtains the following result \cite{Brefine}
\begin{theorem}[Bourgain]
The Shr\"odinger IVP \eqref{shr} in  $\R^{2}$ is globally well-posed in 
$H^{s}$ for $s>3/5$.
\end{theorem}
Using this method several authors extended global well-posedness 
results for variety of  equations, see for example 
\cite{CST99} for the KdV equation, \cite{FLP} for the 
modified KdV, \cite{KPV3} for wave equations, \cite{T} and \cite{Tz}
for the KP-II equation and \cite{takaoka:dnls-global}
 for the Schr\"odinger equations with derivative nonlinearity.

\section{The almost conservation laws: a first attempt}\label{part1}

We  restrict the description of this method to the KdV initial 
value problem \eqref{kdv}. We  remark at the end on the 
applications to other equations. 

To help the reader in understanding this method we decided to 
reproduce in a coherent way the evolution of thoughts that guided us to 
our recent findings. We start by proving the 
conservation of the $L^{2}$-norm for the solution $u$ of \eqref{kdv}
by integration by parts. We refer to this proof as a proof in
{\em physical space} in contrast with another one that we will give 
later and that will be performed in {\em frequency space}. 
If we multiply the equation in \eqref{kdv} by $u$ we obtain 
\[\frac{1}{2}\partial_t u^2 = - \partial_x 
( u \partial_x^2 u ) + \partial_x
\left( \frac{1}{2} [\partial_x u ]^2 \right) 
-\frac{1}{3} \partial_x u^3\]
and integration over the line, or in the periodic case, over the 
circle, we obtain the desired identity
\begin{equation}
    \frac{d}{dt}\|u\|_{L^{2}}^{2}=0.
    \end{equation}
This type of proof does not involve any analysis of the interaction 
of frequencies, which 
we believe is the key to understand the evolution not just of the 
of the $L^{2}$, but also of the $H^{s}$ norms, for any $s\in \R$.

We recall that in Section \ref{con} we observed that the method of 
conservation laws cannot establish global results for \eqref{kdv}
on the line, 
when the initial data $u_{0}\in H^{s}$ for $s\in (-3/4,0)$. So we assume 
that $s<0$. There are no conservation laws, that we are aware of,
for the $H^{s}$ norm, when $s$ is negative, hence some new idea has to be considered. We 
borrow from Bourgain \cite{Brefine} the splitting process 
into 
low and high frequency, but this time  the splitting is done in a 
smooth way and on the solution $u$ itself, not the initial data. 
This argument has been successfully used by Keel and Tao for the 
1D wave map problem \cite{KeeTao98}. So we consider the multiplier
\begin{equation}\label{specialm}
{{\widehat{Iu}}}( \xi ) = m( \xi ) {\widehat{u}} ( \xi ),
~m( \xi ) = \left\{ \begin{matrix}
1, & |\xi | < N, \\
N^{-s} {{|\xi |}^s} , & |\xi | \geq 10 N
\end{matrix}
\right.
\end{equation}
where $m$ is smooth and monotone and $N$ is a large number to be fixed 
later. The operator $I$ (barely) maps $H^s ( \R ) 
\longmapsto L^2 ( \R )$. Observe
that on low frequencies $\{ \xi : |\xi | < N \},~I$ 
is the identity operator.
Note also that $I$ commutes with differential operators. 
We now want to repeat the argument presented above to prove the conservation 
of the $L^{2}$ norm, but this time for $\|Iu(t)\|_{L^{2}}$. 
 Using the Fundamental 
Theorem of Calculus, the equation, and integration by parts, we have
\begin{eqnarray}
\label{l2almost}{{\| Iu ( t) \|}_{L^2}^2 }
&=& {{\| Iu (0 ) \|}_{L^2}^2 } + \int_0^t \frac{d}{d\tau }
( Iu( \tau) , Iu ( \tau ) ) d\tau, \\
\nonumber& = &  {{\| Iu (0 ) \|}_{L^2}^2 } + 2 \int_0^t
( I {\dot {u}} ( \tau ) , Iu( \tau )) d \tau , \\
\nonumber&= & {{\| Iu (0 ) \|}_{L^2}^2 } + 2 \int_0^t
( I ( - {u_{xxx}} - \half \partial_x [u^2] )(\tau) , Iu( \tau )) d\tau \\
\nonumber& = &  {{\| Iu (0 ) \|}_{L^2}^2 } + \int_0^t
(I ( - \partial_x [u^2] ), Iu) d\tau,
\end{eqnarray}
where $(\cdot,\cdot)$ is the scalar product in $L^{2}$.
The error that could make $\|Iu(t)\|_{L^{2}}$ too large in the future
is 
\begin{equation}\label{error}
    R(t)=\int_0^t(I ( - \partial_x [u^2] ), Iu) d\tau.
\end{equation}
The idea is to use local well-posedness estimates to show that
locally in time $R(t)$ is small.  To do so we first have to recall the 
precise local well-posedness result of Kenig, Ponce and Vega 
\cite{KPVBilin}. We define the space $X^{s,b}, \, s,b\in \R$ 
as the closure of the 
Schwartz's functions with respect to the norm 
$$\|f\|_{X^{s,b}}=\left(\int_{\R^{2}}|\hat f|^{2}(\xi,\tau)
(1+|\xi|)^{2s}(1+|\tau-\xi^{3}|)^{2b}\right)^{1/2}.$$
Observe that for $b>1/2$, it follows that 
$X^{s,b}\subset C([0,T],H^{s})$. Kenig, Ponce and Vega proved the 
following bilinear estimate \cite{KPVBilin}
\begin{theorem}[Kenig-Ponce-Vega]
For $s>-3/4$ and $b>1/2$, there exists $b'<b$ such that
\begin{equation}\label{bilinear}
    \|\partial_{x}(uv)\|_{X^{s,b'-1}}\lesssim \|u\|_{X^{s,b}}
\|v\|_{X^{s,b}}.    
\end{equation}
Moreover if $s<-3/4$, there is no $b$ and $b'$ such that 
\eqref{bilinear} is true.
    \end{theorem}
This bilinear estimate is essential to obtain an estimate like 
\eqref{nonlin} and hence to use a fixed point theorem. The 
local well-posedness result can be summarized  in the following theorem. 
Assume that $\psi(t)$ is a cut-off function relative to the interval 
$[-2,2]$, \cite{KPVBilin}.
\begin{theorem}[Kenig-Ponce-Vega]
For any $u_{0}\in H^{s}, s>-3/4$ there exist 
$T=C(\|u_{0}\|_{H^{s}})^{-\alpha}$
and a unique solution $u$ for \eqref{kdv} such that $u$ exists for all 
$t\in [-T,T]$ and in particular
$$
\|\psi(\cdot/T)u\|_{X^{s,b}}\leq C\|u_{0}\|_{H^{s}}.
$$
 \end{theorem}
A modification of this theorem can be proved when we introduce the 
multiplier operator $I$. In fact we have \cite{CKSTTKdV0}
\begin{theorem}\label{localI}
For any $u_{0}\in H^{s}, s>-3/2$ there exist 
$T=C(\|Iu_{0}\|_{L^{2}})^{-\alpha}$
and a unique solution $u$ for \eqref{kdv} such that $u$ exists for all 
$t\in [-T,T]$ and in particular
\begin{equation}\label{local}
\|\psi(\cdot/T)Iu\|_{X^{0,b}}\leq C\|Iu_{0}\|_{L^{2}}.
\end{equation}
 \end{theorem}
Now let's go back to the estimate of the error $R(t)$. Using 
Plancherel  
\begin{eqnarray*}|R(t)|&\leq& 
    \int_{\R}\int_{\R}|\widehat{\partial_{x}I(u^{2})}|(\xi,\tau)
    (1+|\xi|)^{s}(1+|\tau-\xi^{3}|)^{b'-1}\\
    &\times&
    |\widehat{\chi_{t}Iu}|(\xi,\tau)(1+|\xi|)^{-s}
    (1+|\tau-\xi^{3}|)^{-b'+1}\,d\xi\,d\tau,
\end{eqnarray*}    
where $\chi_{t}$ is the characteristic function of $[0,t]$. Then    
 by the Cauchy-Schwarz inequality we have
\begin{equation}\label{error1}
|R(t)|\leq \| \partial_{x}I(u^{2})\|_{X^{0,1-b'}}\|Iu\|_{X^{0,1-b'}}.   
    \end{equation}
Hence if we could prove a bilinear inequality like
\begin{equation}\label{badbil}
    \| \partial_{x}I(u^{2})\|_{X^{0,1-b'}}\leq N^{-\beta}
\|Iu\|_{X^{0,b}}^{2},
\end{equation}
for some $\beta>0$, then we would be done because the factor 
$N^{-\beta}$, would make the error small
\footnote{This will be explained in more details in 
Theorem \ref{gwp}.}. But unfortunately, even 
though  \eqref{badbil} looks a lot like \eqref{bilinear},
 it is false
 \footnote{This is not obvious at first sight, for more explanation one 
 should consult \cite{CKSTTKdV0}.}
 due to the interaction of very low 
frequencies ($|\xi|<<N$) with very large frequencies 
($|\xi|>>N$). But not everything is lost, in fact we can introduce 
for free a 
suitable cancellation\footnote{This cancellation is recognaseble once 
one writes the expression in Fourier transform and uses the mean value 
theorem, see \cite{CKSTTKdV0} for a precise calculation.} by  rewriting
\eqref{error} as 
\begin{equation}\label{error2}
    R(t)=\int_0^t \int\partial_x \left\{
{{(I(u))}^2} - I( u^2 ) \right\}~ Iu ~dx  d\tau,
\end{equation}
and we replace \eqref{error1} with
\begin{equation}\label{error3}
|R(t)|\leq \|\partial_x \{ {{(I(u))}^2} - I( u^2 ) \}
\|_{X^{0,1-b'}}\|Iu\|_{X^{0,1-b'}}.   
    \end{equation}
Now the following desired proposition is true (see \cite{CKSTTKdV0})
\begin{proposition} \label{extrasmoo}(Extra smoothing)
The bilinear estimate
\begin{equation}
\label{smoothingest}
\|\partial_x \{ I(u) I(v) - I( uv )\|_{X^{0,-1/2-}}
\leq C N^{-\frac{3}{4}+}\|Iu\|_{X^{0,1/2+}}
\|Iv\|_{X^{0,1/2+}}.
\end{equation}
holds.
\end{proposition} 
Combining \eqref{l2almost} with \eqref{error3} and  
\eqref{smoothingest}, we obtain the almost conservation 
law\footnote{We refer to these types of estimate as almost conservation 
laws because of the presence of the decaying factor $N^{-\beta}$.}
\begin{equation}\label{l2almostI}
   {{\| Iu ( t) \|}_{L^2}^2 }\leq {{\| Iu ( 0) \|}_{L^2}^2 }
   +CN^{-\frac{3}{4}+}\|Iu\|_{X^{0,1/2+}}^{3}.
    \end{equation}
In proving the following theorem we  describe in detail how one 
obtains a global result by an iteration based on \eqref{l2almostI}.
\begin{theorem}\label{gwp}
 The initial value problem  \eqref{kdv} is globally well posed in 
 $H^{s}$ for all $s$ such that $s>-3/10$.
    \end{theorem}
\begin{proof}
The proof is taken from \cite{CKSTTKdV0}.
Global well-posedness of \eqref{kdv} will
follow if we show well-posedness on $[0,T]$ for arbitrary $T>0$. We 
renormalize things a bit via scaling. If $u$ solves \eqref{kdv} then
$u_\lam ( x ,t ) = {{(\frac{1}{\lam}  )}^2} 
u( \frac{x}{\lam} , \frac{t}{\lam^3})  $ solves \eqref{kdv} with 
initial data 
\begin{equation}\label{scal0}
    u_{0,\lam} (x,t) = 
 {{\left(\frac{1}{\lam}  \right)}^2} u_{0} \left( \frac{x}{\lam} 
 \right).
 \end{equation}
Note that $u$ exists on $[0,T]$ if and only if $u_\lam$ exists on 
$[0,\lam^3 T]$. A calculation shows that
\begin{equation}
\|Iu_{0,\lam} \|_{L^{2}} \leq C {\lam^{-\frac{3}{2} - s }} N^{-s}  
\|u_{0}\|_{H^{s}} .
\label{scaling}
\end{equation}
Here $N= N(T)$ will be selected later but we choose $\lam = 
\lam(N)$ right now by requiring
\begin{equation}
 C {\lam^{-\frac{3}{2} - s }} N^{-s}  \|u_{0}\|_{H^{s}}
 \thicksim 1 \implies
\lam \thicksim N^{- \frac{2s}{3 + 2s }}.
\label{lamchosen}
\end{equation}
We now drop the $\lam$ subscript on $u_{0}$ by assuming that 
\begin{equation}
\|I u_{0}\|_{L^{2}}= \epsilon_0 \ll 1,
\label{normalized}
\end{equation}
and our goal is to construct the solution of \eqref{kdv} on the time
interval $[0, \lam^3 T]$.

The local well-posedness result of Theorem \ref{localI}  shows we 
can construct the
solution for $t \in [0, 1]$ if we choose $\epsilon_0$ small enough. 
Using 
\eqref{local} and \eqref{normalized}, the almost $L^2$ 
conservation property \eqref{l2almostI} we obtain
\begin{equation*}
{{\| Iu(1) \|}_2^2} \leq \e_0^2 + N^{-\frac{3}{4}+}.
\end{equation*}
We can iterate this process $N^{\frac{3}{4}-}$ times before doubling 
$\|Iu(t)\|_{L^{2}}$. Therefore, we advance the solution by taking 
$N^{\frac{3}{4}-}$ time steps of size $O(1)$. We now restrict 
$s$ by demanding that
\begin{equation}
N^{\frac{3}{4}-} \gtrsim \lam^3 T = N^{\frac{-6s}{3 + 2 s}} T
\label{scondition}
\end{equation}
is ensured for large enough $N$,
so $s> - \frac{3}{10}$.
\end{proof}

\section{The almost conservation laws: the final version}\label{part2}

The cancellation that we introduced in \eqref{error3}, and that
can be seen explicitly in frequency space by taking Fourier 
transforms, led us to try to understand more deeply the interaction of 
frequencies during the evolution of the solution $u(x,t)$ of 
\eqref{kdv}. For this purpose we 
propose here a proof in frequency space of the $L^{2}$ conservation law 
for the solution of \eqref{kdv}. By the 
Plancherel theorem we have
\begin{eqnarray*}
& &\|u(t)\|_{L^{2}}^{2}=\int {\widehat{u}} ( \xi ) 
{\overline{\widehat{u}}} 
( \xi ) d\xi
=\int  {\widehat{u}} ( \xi ) {\widehat{\overline{u}}} ( -\xi ) d\xi
= \int_{\xi_{1}+\xi_{2}=0} {\widehat{u}} ( \xi_1 ) 
{\widehat{u}} ( \xi_2 )\,d\xi_{1}\,d\xi_{2},
\end{eqnarray*}
since $u$ is $\R$-valued. Therefore, by substituting in the equation
we obtain
\begin{eqnarray*}\partial_t \int {\widehat{u}} ( \xi ) 
{\overline{\widehat{u}}} ( \xi ) d\xi 
&=& 2 \int_{\xi_{1}+\xi_{2}=0} {\widehat{u_t}} ( \xi_1 ) {\widehat{u}} 
(\xi_2 )\,d\xi_{1}\,d\xi_{2}\\
&=& 2 \int_{\xi_{1}+\xi_{2}=0} \left[ - (i \xi_1 )^3 
{\widehat{u}} (\xi_1)
- \frac{1}{2} (i \xi_1 ) {\widehat{u^2}} ( \xi_1) 
\right] {\widehat{u}} ( 
\xi_2 )\,d\xi_{1}\,d\xi_{2}.
\end{eqnarray*}
Now we symmetrize the first term and we expand the convolution to get
\begin{eqnarray*}
\frac{d}{dt}\int u^{2}(x)\,dx&=&\partial_t \int {\widehat{u}} ( \xi ) 
 {\overline{\widehat{u}}} ( \xi ) d\xi
 =-\int_{\xi_{1}+\xi_{2}=0} i ( \xi_1^3 + \xi_2^3 ) {\widehat{u}} 
 ( \xi_1 )
{\widehat{u}} (\xi_2 )\,d\xi_{1}\,d\xi_{2}\\
&&- \int_{\xi_{1}+\xi_{2}+\xi_{3}=0} i 
(\xi_1 + \xi_2 )
{\widehat{u}} (\xi_1 ){\widehat{u}} (\xi_2 ){\widehat{u}} (\xi_3 )
\,d\xi_{1}\,d\xi_{2}\,d\xi_{3}.
\end{eqnarray*}
The first term is clearly zero. Upon writing $\xi_1 + \xi_2 = -\xi_3$
and symmetrizing, the second term vanishes too.
This  {\em symmetrization/cancellation} describes the non linear 
interaction 
of the  frequencies of the solution $u$ for the KdV equation
in \eqref{kdv}.  We stress here once more that we think this is 
an important mechanism to understand in order to keep track of the 
various pieces of $\hat{u}$ once we perform a frequency 
localization like we did by introducing the multiplier operator $I$.

It is time now to introduce some notation that will make the rest of 
our presentation 
less cumbersome. We start with the following definitions:

\begin{definition} 
A $k$-multiplier $m$ is a function
$m: \R^k \longrightarrow \C$. A $k$-multiplier 
is symmetric if
$m( \xi ) = m( \sigma( \xi))$ for all $\sigma \in S_k $. The 
symmetrization of a $k$-multiplier is
\[[m]_{sym} ( \xi ) = \frac{1}{n!} \sum\limits_{\sigma \in S_k }
m(\sigma(\xi)).\]
A $k$-multiplier  generates the $k$-linear functional  
via the integration
\begin{equation}
    \Lambda_k (m ) = \int_{A_{k}}
m( \xi_1 , \dots \xi_k ) {\widehat{u}} (\xi_1 ) \dots {\widehat{u}} 
(\xi_k ),
\end{equation}
where $A_{k}=\{(\xi_{1},\ldots,\xi_{k}) / \xi_1 + \dots +\xi_k = 0\}$.
\end{definition}
We immediately observe that we can rewrite $\|Iu\|_{L^{2}}$ using 
the $\Lambda$ notation above. In fact 
$$\|Iu(t)\|_{L^{2}}^{2}=\int_{A_{2}}m(\xi_{1}m(\xi_{2})\hat u(\xi_{1},t)
\hat u(\xi_{2},t)=\Lambda_{2}(m(\xi_{1})m(\xi_{2})).$$
It is then clear the purpose of next proposition:

\begin{proposition}\label{der}
Suppose $u$ satisfies the KdV equation, and $m$ is a symmetric 
$k$-multiplier and 
\[\Lambda_k (m ) = \int_{A_{k}}
m( \xi_1 , \dots \xi_k ) {\widehat{u}} (\xi_1 ) \dots {\widehat{u}} 
(\xi_k ),\]
is the $k$-linear functional generated by $m$. Then
\begin{equation}
\frac{d}{dt} \Lambda_k (m )=\Lambda_k (m \alpha_k ) 
- i \frac{k}{2}
\Lambda_{k+1} \left(\tilde{m}(\xi_1, \dots,\xi_{k+1})\right),
\end{equation}
where 
\[\alpha_k = i (\xi_1^3 + \dots + \xi_k^3 )\]
and
$$
\tilde{m}(\xi_{1},\ldots,\xi_{k+1})=
m(\xi_1,\dots,\xi_{k-1}, [\xi_k + \xi_{k+1}])(\xi_k + \xi_{k+1}).
$$
\end{proposition}
We  now describe the general principle 
behind the almost conservation laws.

Let $m$ be an $\R$-valued even 1-multiplier. Define again 
the multiplier operator\footnote{This is the same operator introduced 
in Section \ref{part1} when we take $m$ as in \eqref{specialm}.} 
$I$ via
\[{\widehat{If}} ( \xi) = m(\xi) {\widehat{f}} (\xi ).\]
For convenience of notation we rename 
$E^{1}_{I}(t)=\|u(t)\|_{H^{s}}^{2}$ and 
\[E_I^2 ( t) = {{\| Iu (t) \|}_{L^2 }^2 } = 
\Lambda_2 ( m(\xi_1 )m(\xi_2)).\]
Our goal now is to define a hierarchy of 
 modified  energies $E_I^i(t), i=2,3,\ldots$ for the solution of 
the IVP \eqref{kdv} such that, when $m$ is like in \eqref{specialm},
\begin{eqnarray}
\label{compareson}& &\|u(t)\|_{H^{s}}\lesssim
E_I^i(t)\lesssim E_I^{i+1}(t),\\
 \label{impro}& &(E_I^{i+1}(b)-E_I^{i+1}(a))<<(E_I^i(b)-E_I^i(a)),
\end{eqnarray}
for any fixed interval $[a,b]$.
In other words we want to find better generations of 
energies that are comparable to the original norm 
$\|u(t)\|_{H^{s}}$, but which increments decrease as the generations 
evolve.

We now present an algorithm that formally\footnote{At this stage not 
all the mathematical quantities that we write 
are proven to make sense, we will worry about this later, for the 
moment we just want to give the original flow of ideas that brought us 
to  rigorous and useful results!} provides improved generations
of energies. 
Using Proposition \ref{der} we calculate
\begin{eqnarray*}
\frac{d}{dt}E^2_I (t)&=&\Lambda_2 ( m(\xi_1 ) m(\xi_2 ) \alpha_2 )
- i \Lambda_3 ( m(\xi_1 ) m(\xi_2 + \xi_3 ) [\xi_2 + \xi_3 ])\\
&=&\Lambda_3 ( -i [ m(\xi_1 ) m(\xi_2 + \xi_3 ) (\xi_2 + \xi_3 ) 
]_{sym} ).
\end{eqnarray*}
We should point out that for $m$ as in \eqref{specialm}
$$R(t)=\int_{0}^{t}
\Lambda_3 ( -i [ m(\xi_1 ) m(\xi_2 + \xi_3 ) (\xi_2 + \xi_3 ) 
]_{sym} )\,ds$$
where $R(t)$ is the error defined in Section \ref{part1}.
We proved 
in Proposition \ref{extrasmoo} that even though $R(t)$ is a 
threelinear expression coming from a bilinear expression such as 
$E^{2}_{I}$, the symmetrization\footnote{This is what we 
called then {\em cancellation}.} allows us to obtain a decay in $N$ 
which is in fact what gives \eqref{impro}. So our goal is to push this 
idea further in the following way:
we first denote
\begin{equation}\label{m3}
    M_3 ( \xi_1 , \xi_2 , \xi_3 )=-i ([ m(\xi_1 ) m(\xi_2 + \xi_3 ) 
(\xi_2 + \xi_3 ) ]_{sym} ).
\end{equation}
Then we define the third generation of modified energy as
\[E_I^3 (t) = E_I^2 ( t) + \Lambda_3 (\sigma_3 ),\]
where $\sigma_{3}$ is a multiplier that will be chosen later.
Now again by Proposition \ref{der} we have
\[\frac{d}{dt}
E_I^3 (t) = \Lambda_3 (M_3 ) + \Lambda_{3}( \sigma_3 \alpha_3 ) +
\Lambda_4 (M_4),\]
where 
\begin{equation}\label{m4}
    M_{4}(\xi_{1},\ldots,\xi_{4})=\sigma_{3}(\xi_{3}+\xi_{4}).
    \end{equation}
We choose $\sigma_{3}$ to cancel the $\Lambda_3$ terms, that is 
\[\sigma_3 = - \frac{M_3}{\alpha_3}.\]
Because $\alpha_{3}=\sum_{i=1}^{3}\xi_{i}^{3}$, 
we expect that
\begin{equation}
    \label{hope}|M_{4}|=\left|M_3\frac{(\xi_{3}+\xi_{4})}{\alpha_3}
    \right|\ll |M_{3}|,
    \end{equation}
    and hence \eqref{impro}.
Certainly at this point our expectation is a pure leap of faith 
because anybody could argue that when\footnote{
 The reader should observe that $\sum_{i=1}^{3}\xi_{i}^{3}=0$
 is the relationship that defines the resonance set of three wave 
 interaction!} $\sum_{i=1}^{3}\xi_{i}^{3}=0$,
the left hand side of the expression in \eqref{hope} would become infinity unless a miraculous
cancellation occurs in the numerator. What really amazed us was that 
indeed such a miracle happens! The ``miracle'' is a combination of the 
type of 
frequency cancellation that we observed in the proof of the $L^{2}$ 
conservation law via the frequency method,  with several applications 
of the Mean Value Theorem that we can perform since we are assuming 
that the multiplier $m$ is smooth, see \cite{CKSTTKdV} for details.

The process we described above may be iterated to formally generate 
a sequence of modified energies
$\{ E^j_I (t) \}_{j=2}^\infty $, with the property that
\[\frac{d}{dt}
E^j_I (t) = \Lambda_{j+1} (M_{j+1} ).\]
The hard part of the argument is to present a rigorous 
proof for the statement
\begin{equation}\label{miracle}
|M_{j+1} |<<|M_{j} |\end{equation}
in an appropriate sense! 

Before we proceed to a less formal, but more technical discussion on 
the algorithm  above, we want to convince the reader that 
in principle our 
method could be used to recover all the conservation laws that the KdV 
equation enjoys. We didn't set to the onerous task of checking this 
in detail, but we can show at least an example that is not trivial,
see also \cite{CKSTTKdV}, the paper where this computation first 
appeared.

We first specify the multiplier $m$ by setting 
$m(\xi)=i\xi$. Then 
\begin{equation*}
  E^2(t) = {{\| \partial_x u \|}_{L^2}^2} = \Lam_2 ( (i\xi_1 ) 
(i \xi_2)).
\end{equation*}
Next we define $E^3(t) = E^2 (t) + \Lam_3 ( \sigma_3 )$, 
and we use Proposition \ref{der} to see that
\begin{equation*}
  \partial_t E^3 (t) = \Lam_3 ( [ i(\xi_1 + \xi_2 ) i \xi_3 \{ 
\xi_1 + \xi_2 \}
]_{sym} ) + \Lam_3 ( \sigma_3 \alpha_3 ) + \Lambda_4 (M_4),
\end{equation*}
where $M_4$ is explicitly obtained from $\sigma_3$. Noting that
$i (\xi_1 + \xi_2 ) i \xi_3 \{ \xi_1 + \xi_2 \} = - \xi_3^3 $ on the 
set 
$\xi_1 + \xi_2 + \xi_3 =0$, we know that
\begin{equation*}
  \partial_t E^3 
(t) = \Lam_3 (-\frac{1}{3} \alpha_3 ) + \Lam_3 (\sigma_3 \alpha_3 ) + 
\Lam_4 ( M_4 ).
\end{equation*}
The choice of $\sigma_3 = \frac{1}{3}$ results in a cancellation of 
the
$\Lam_3$ terms and
$$M_4 = [\{ \xi_1 + \xi_2 \}]_{sym} = \xi_1 + \xi_2 
+ \xi_3 
+ \xi_4,$$
so $M_4 = 0$. Therefore, $E^3 
(t) = \Lam_2 ( (i \xi_1 ) (i \xi_2 )) + \Lam_3 ( \frac{1}{3} )$
is an exactly conserved quantity. The modified energy construction 
applied
to the Dirichlet energy ${{\| \partial_x u \|}_{L^2}^2}$ led us to the Hamiltonian for KdV
described in \eqref{hamkdv}. Applying the
construction to higher order derivatives in $L^2$ we expect that 
will similarly lead to the higher conservation laws of KdV.

Assume now that the initial data $u_{0}$ of our IVP is in $H^{s}, \, s
\in (-3/4,0)$. Let $m$ be the multiplier defined in \eqref{specialm}.

Using multilinear type estimate one can  show that
\cite{CKSTTKdV}
\begin{equation}\label{gn}
    \|u(t)\|_{H^{s}}\lesssim E^2_I (t) \lesssim E^4_I (t).
    \end{equation}
But the hart of the matter is the following proposition
\begin{proposition} For fixed $T>0$
\begin{equation}\label{sharp}
E^4_I ( T ) - E^4_I (0) = \int_0^T \Lambda_5 ( M_5 )(\tau ) d 
\tau\leq C_{T} N^{-3  + \epsilon} {{\| Iu \|}_{X^{0,1/2+}}^5}.
\end{equation}
\end{proposition}
For a complete proof see \cite{CKSTTKdV}.

At this point probably the reader would like to ask the following 
question: Why did we stop at $E^{4}_{I}$? The obvious answer that we 
can give  is that we stopped 
because the decay of the increment of this modified energy, given by 
\eqref{sharp}, is enough to obtain  the best possible  result:
\begin{theorem}\label{bestkdv}
 The IVP \eqref{kdv}
is globally well-posed  in $H^{s}$ for $s>-3/4$.   
    \end{theorem}
But  there is a much deeper reason why  we didn't pursue the 
estimates of the increment of the energies $E^{k}_{I}$, for $k>4$.   The 
formal expression for the increment of these 
energies  becomes more and more 
complex. Nice algebraic properties like \eqref{factthree} and 
\eqref{factfour} below
are no longer available!
Also it seems  to us that the reason why we didn't need to estimate 
the increment for all the modified energies is that $-3/4$ is larger 
than the scaling index\footnote{The scaling index is the 
the Sobolev index $s_{c}$ such that the rescaled initial data 
$u_{0,\lam}$ defined in \eqref{scal0} has the property that
$\|u_{0,\lam}\|_{\dot{H}^{s_{c}}}$ is independent of $\lam$.}, 
which, in this case, is $-3/2$. 

The proof of Theorem \ref{bestkdv} is similar to the proof of 
Theorem  \ref{gwp} if one uses \eqref{gn} and replaces 
\eqref{l2almostI} with 
\eqref{sharp}, see \cite{CKSTTKdV} for details.

To give an idea of the type 
of miracle that makes \eqref{miracle} analytically correct
we  consider $M_{4}$, defined in \eqref{m4}. The complete 
estimate of $M_{4}$ is very involved, so we will restrict 
ourselves to some special cases. 
The computations that follow are directly taken from 
\cite{CKSTTKdV}.
 We  recall the following
arithmetic facts that may be easily verified:
\begin{equation}
\xi_1 + \xi_2 + \xi_3 = 0 \implies \alpha_3 = \xi_1^3 + \xi_2^3
+ \xi_3^3 = 3 \xi_1 \xi_2 \xi_3.
\label{factthree}
\end{equation}
\begin{equation}
\xi_1 + \xi_2 + \xi_3 + \xi_4 = 0 \implies \alpha_4 = \xi_1^3 + 
\xi_2^3
+ \xi_3^3 + \xi_4^3 = 3 (\xi_1 + \xi_2) ( \xi_1 + \xi_3 ) (\xi_1 + 
\xi_4).
\label{factfour}
\end{equation}
Recall that, 
\begin{equation}
  \label{Mfour}
  M_4 ( \xi_1 , \xi_2 , \xi_3 , \xi_4 ) = c [ \sigma_3 ( \xi_1 , 
\xi_2, \xi_3
+ \xi_4 ) (\xi_3 + \xi_4 ) ]_{sym},
\end{equation}
where $\sigma_3 = - \frac{M_3}{\alpha_3}$ and
\begin{eqnarray}
  \label{Mthree}
  M_3 ( x_1 , x_2 , x_3 ) &=& -i [ m(x_1 ) m( x_2 + x_3 )(x_2 + x_3 ) 
]_{sym}\\
\nonumber&=& - \frac{i}{3}[m^2 (x_1 ) x_1 + m^2 (x_2) x_2 + m^2 (x_3 ) x_3] , 
\end{eqnarray}
and by \eqref{factthree}
$\alpha_3 (x_1 , x_2 , x_3 ) = x_1^3 + x_2^3 + x_3^3 = 3 x_1 x_2 
x_3 $.
We shall ignore the irrelevant constant in \eqref{Mfour}.
Therefore,
\begin{equation}
\label{Mfromm}
   M_4 ( \xi_1 , \xi_2 , \xi_3 , \xi_4 ) = -\half \left[ \frac{ 
m^2 (\xi_1 ) \xi_1  + 
m^2 (\xi_2) \xi_2 + m^2 (\xi_3 + \xi_4 )(\xi_3 + \xi_4 ) }{
3 \xi_1 \xi_2 } \right]_{sym}  
\end{equation}
\begin{equation*}
= -\half \left[ \frac{ 
2 m^2 (\xi_1 ) \xi_1  + m^2 (\xi_3 + \xi_4 )(\xi_3 + \xi_4 ) }{
3 \xi_1 \xi_2 } \right]_{sym} .
\end{equation*}
Using the identity \eqref{factfour} and lots of symmetrizations and 
clever tricks like in  \cite{CKSTTKdV}, one can reexpress $M_{4}$ as 
\begin{equation}
  \label{oneten}
  M_4 ( \xi_1 , \xi_2 , \xi_3 , \xi_4 )
= -\frac{1}{36} \frac{1}{\xi_1 \xi_2 \xi_3 \xi_4 } \times
\end{equation}
\begin{eqnarray*}
& \{ & \xi_1 \xi_2 \xi_3 [ m^2 ( \xi_1 ) + m^2 ( \xi_2 ) + m^2 ( 
\xi_3 ) 
- m^2 ( \xi_1 + \xi_2 ) - m^2 ( \xi_1 + \xi_3 ) - m^2 ( \xi_1 + \xi_4 
) ]  \\
&+& \xi_1 \xi_2 \xi_4 [ m^2 ( \xi_1 ) + m^2 ( \xi_2 ) + m^2 ( \xi_4 ) 
- m^2 ( \xi_1 + \xi_2 ) - m^2 ( \xi_1 + \xi_3 ) - m^2 ( \xi_1 + \xi_4 
) ] \\
&+&
 \xi_1 \xi_3 \xi_4 [ m^2 ( \xi_1 ) + m^2 ( \xi_3 ) + m^2 ( \xi_4 ) 
- m^2 ( \xi_1 + \xi_2 ) - m^2 ( \xi_1 + \xi_3 ) - m^2 ( \xi_1 + \xi_4 
) ] \\
&+& \xi_2 \xi_3 \xi_4 [ m^2 ( \xi_2 ) + m^2 ( \xi_3 ) + m^2 ( \xi_4 ) 
- m^2 ( \xi_1 + \xi_2 ) - m^2 ( \xi_1 + \xi_3 ) - m^2 ( \xi_1 + \xi_4 
) ] \}.\\
\end{eqnarray*}
Assume now that $m$ is like in \eqref{specialm} and that 
$\xi_{i}=0$ for $i=1,\ldots,4$. Then obviously $M_{4}=0$. To make 
things more interesting let's now assume that  only $\xi_{1}=0$.
Then the numerator of $M_{4}$ takes the form of
$$\xi_2 \xi_3 \xi_4 
[ m^2 ( \xi_2 ) + m^2 ( \xi_3 ) + m^2 ( \xi_4 ) 
- m^2 ( \xi_2 ) - m^2 ( \xi_3 ) - m^2 ( \xi_4 ) ] \}$$
which is once again zero.

We end this section and the article with some general remarks.
Using the arguments presented in this section we are able to 
completely fill the gap between local well-posedness  and 
global well-posedness also for the periodic KdV and 
the continuous and periodic  mKdV. The periodic KdV
problem is  more difficult because the scaling argument used 
in the proof of Theorem \ref{bestkdv} changes the period of the 
 rescaled solution, hence all the estimates have to be
independent of the rescaling parameter $\lambda$ up to 
a factor $\lambda^{0+}$. To 
approach the mKdV problem we use the Miura transformation that 
relates solutions of the KdV to solution of the mKdV equation 
in an explicit way. For details the reader should see \cite{CKSTTKdV}.

The method of almost conservation laws that we presented  here 
is very general. We used it to obtain similar sharp results for the 
1D Schr\"odinger equation with derivative nonlinearity 
\cite{CKSTTDNLS} \cite{CKSTTDNLS1}, and to obtain 
partial results for the IVP \eqref{shr}
\cite{CKSTTNLS}, that improve Bourgain's results in \cite{Brefine}. 

We believe that given a dispersive equation, the method 
we developed gives an analytic tool to study the nonlinear 
interactions of  parts of the solution of the equation 
carried by different frequencies. 
We are now entering the domain  of the weak turbulence theory!

\bibliographystyle{amsalpha}

\enddocument